\documentclass[11pt]{article}
\usepackage{amsfonts}
\usepackage{amsmath}
\usepackage{amssymb}
\usepackage{cite}

\newtheorem{theo}{Theorem}
\newtheorem{exam}[theo]{Example}
\newtheorem{lem} [theo]{Lemma}

\makeatletter \@addtoreset{equation}{section}
\@addtoreset{theo}{}\makeatother

\makeindex 
\def\pf{\noindent {\it Proof.} }
\def\qed{\hfill \rule{4pt}{7pt}}
\def\n{\boldsymbol}
\def\Ann{\mathop{\rm Ann}}
\def\res{\mathop{\rm res}}
\def\bfal{{\boldsymbol{\alpha}}}
\begin{document}

\title{\textbf{Formal residue and computer proofs of combinatorial identities }}
\author{Qing-Hu Hou$^1$, and Hai-Tao Jin$^2$\\
\small Center for Combinatorics, LPMC-TJKLC\\[-0.8ex]
\small Nankai University, Tianjin 300071, P.R. China\\
\small \texttt{$^1$hou@nankai.edu.cn, $^2$jinht1006@mail.nankai.edu.cn}
}
\date{}
\maketitle

\begin{abstract}
The coefficient of $x^{-1}$ of a formal Laurent series $f(x)$ is called the formal residue of $f(x)$.
Many combinatorial numbers can be represented by the formal residues of hypergeometric terms. With these representations and the extended Zeilberger's algorithm, we generate recurrence relations for summations involving combinatorial sequences such as Stirling numbers and their $q$-analogue. As examples, we give computer proofs of several known identities and derive some new identities. The applicability of this method is also studied.
\end{abstract}



\section{Introduction}
Finding recurrence relations for summations is a key step in computer proofs of combinatorial identities. In 1990's, Wilf and Zeilberger \cite{WilfZeilberger1992,Zeilberger1991} developed
the method of creative telescoping to generate recurrence relations for hypergeometric summations. Since then, many extensions and new algorithms
have been discovered and designed for various kinds of summations.
See, for example, \cite{Chyzak1998,Chyzak2000} for holonomic sequences, \cite{Wegschaider1997,ChenHouMu2010} for
multivariable hypergeometric terms, \cite{Schneider2001} for nested sums and products, \cite{Kauers2007,KauersSchneider} for Stirling-like numbers, \cite{Chyzak-Kauers-Balvy2009,Koutschan2009,Mu2011}
for non-holonomic sequences.

Our approach is motivated by the work of Chen and Sun \cite{ChenSun2009}. By using
the Cauchy contour integral representations, they transformed sums involving Bernoulli numbers into hypergeometric summations. Then the recurrence relations for the sums can be derived by the extended Zeilberger's algorithm \cite{ChenHouMu2010}.

In the present paper, we combine the formal residue operator and the extended Zeilberger's algorithm to generate recurrence relations for combinatorial sums. With this residue method, we give computer proofs of some known identities and derive some new identities.
Moreover, we study the applicability of this method. We show that in the case of one variable, it is equivalent to the Sister-Celine's method.

We note that Egorychev \cite{Egorychev1984} provided integral representations for many combinatorial numbers and used them to prove combinatorial identities. F\"urst \cite{Christoph2011} reformulated Egorychev's method in terms of formal residue operators. Egorychev transformed sums into geometric series and then evaluated them by some manipulation rules. We transform sums into hypergeometric sums and find recurrence relations they satisfied.

The paper is organized as follows. In Section 2, we describe the residue method. Then in Section 3, we give
several examples involving Stirling numbers of both kinds. Section 4 is devoted to deriving
two new Stirling number identities. In section 5, we consider
the $q$-Stirling numbers as well as other combinatorial sequences which also fall in the scope of
our method. Finally in Section 6, we study the applicability of the residue method.


\section{The method of residue}

Let $\mathbb{K}$ be a field and $\mathbb{K}((z))$ be the set of formal Laurent series in
the indeterminate $z$ over $\mathbb{K}$. For any element
\begin{equation}\label{fls}
f(z) = \sum_{n=n_0}^\infty a_n z^n \in \mathbb{K}((z)),
\end{equation}
the {\it formal residue operator} $\res\limits_z$ (or, $\res$ if no confusion) is defined by
\[
\res f(z) = \res_z f(z) = a_{-1}.
\]
Clearly, the $k$-th coefficient of $f(z)$ can be represented by formal residue as follows
\[
a_k=\res \frac{f(z)}{z^{k+1}}.
\]
We see that this representation is equivalent to the Cauchy integral representation of $a_k$,
\[
a_k=\frac{1}{2\pi i}\oint _{|z|=\rho}\frac{f(z)}{z^{k+1}}dz.
\]

Based on the formal residue, we give a computer assisted method to derive recurrence relations for
sums involving non-hypergeometric sequences. Consider a definite sum with the form of
\[
f(\n{n}) = \sum_{k=-\infty}^\infty F(\n{n}, k),
\]
where $\n{n}=(n_1,\ldots,n_r)$ is the vector of parameters. The residue method consists of the following three steps.
\begin{enumerate}
\item[1.] Rewrite the summand $F(\n{n},k)$ as $\res\limits_z \tilde{F}(\n{n},k,z)$, where $\tilde{F}(\n{n},k,z)$ is a hypergeometric term.

\item[2.] Take a finite subset $S \subset \mathbb{N}^r$ and apply the extended Zeilberger's algorithm to the similar terms $\{ \tilde{F}(\n{n}+\bfal,k,z) \}_{\bfal \in S}$, where $\mathbb{N}$ denotes the set of non-negative integers. If the algorithm succeeds, we thus obtain a relation of the form
\begin{equation}\label{tF-rec}
\sum_{\bfal \in S} p_\bfal(\n{n}) \tilde{F}(\n{n}+\bfal, k, z) = \Delta_k G(\n{n},k,z),
\end{equation}
where $p_\bfal(\n{n})$ are polynomial coefficients independent of $k$ and $z$ and $G(\n{n},k,z)$ is a hypergeometric term similar to $\tilde{F}(\n{n}, k, z)$.

\item[3.] Summing over $k$ and applying the operator $\res_z$, we are led to a recurrence relation for the sum $f(\n{n})$,
\[
\sum_{\bfal \in S} p_\bfal(\n{n}) f(\n{n}+\bfal) = \res_z G(\n{n},+\infty,z) - \res_z G(\n{n},-\infty,z).
\]
\end{enumerate}
\noindent \emph{Remark.} In most cases, $G(\n{n},k,z)$ is finite supported and hence we do not need
to calculate $\res_z G(\n{n},+\infty,z)$ and $\res_z G(\n{n},-\infty,z)$.

To conclude this section, we give an example to illustrate the method of residue. More examples can be found in Sections 3--5.

\begin{exam}We have \cite[identity (6.15)]{GKP1994}
\begin{equation}\label{GKP6.15}
\sum_{k}\binom{n}{k}S_2(k,m)=S_2(n+1,m+1),
\end{equation}
where $S_2(n,m)$ is the Stirling number of the second kind.
\end{exam}
\pf
It is known that
\[
\sum_{n=k}^{\infty}S_2(n,k)z^n=\frac{z^k}{(1-z)(1-2z)\cdots(1-kz)}.
\]
Therefore,
\begin{equation}\label{Stiring second res}
S_2(n,k)=\res_z \frac{z^k}{z^{n+1}(1-z)(1-2z)\cdots(1-k z)}.
\end{equation}
Denote the left hand side of \eqref{GKP6.15} by $L(n,m)$. We thus have
\[
L(n,m)=\res_z\sum_{k}\binom{n}{k}\frac{z^m}{z^{k+1}\prod_{i=1}^{m}(1-iz)}.
\]

Now consider the inner summand
\[
C(n,m,k,z)=\binom{n}{k}\frac{z^m}{z^{k+1}\prod_{i=1}^{m}(1-iz)}.
\]
Applying the extended Zeilberger's algorithm to the four similar terms
\[
C(n+i, m+j, k, z), \quad 0 \le i,j \le 1,
\]
we find that
\[
C(n+1,m+1,k,z)-(m+2)C(n,m+1,k,z)-C(n,m,k,z)=\Delta_k \frac{-kzC(n,m,k,z)}{(n+1-k)(1-(m+1)z)},
\]
Summing over $k$ and applying the formal residue operator, we derive that
\[
-L(n,m)-(m+2)L(n,m+1)+L(n+1,m+1)=0.
\]
This agrees with the recurrence relation satisfied by $S_2(n+1,m+1)$.
Finally, the identity follows by checking the initial values
\[
L(0,m)=S_2(0,m)=S_2(1,m+1), \quad L(n,0)=1=S_2(n+1,1). \tag*{\qed}
\]

We remark that most of the sums appearing in this paper can also be treated by
Koutschan's implementation of the creative telescoping algorithm on non-holonomic sequence (for more detail,
see \cite{Koutschan2009}).
The only exception is Example 2.

\section{Stirling number identities}\label{Stirling section}
In this section we shall provide several examples involving Stirling numbers of both kinds to
illustrate the residue method. Recall that
\[
\sum_{k=0}^{n}S_1(n,k)z^k=(z)_{\underline{n}} = z(z-1)\cdots(z-n+1),
\]
and
\[
\sum_{n=k}^{\infty}S_2(n,k)z^n=\frac{z^k}{(1-z)(1-2z)\cdots(1-kz)},
\]
where $S_1(n,k)$ and $S_2(n,k)$ are Stirling numbers of the first kind and of the second kind, respectively. We thus have
\[
S_1(n,k)=\res_z \frac{(z)_{\underline{n}}}{z^{k+1}}
\]
and
\[
S_2(n,k)=\res_z \frac{z^k}{z^{n+1}(1-z)(1-2z)\cdots(1-k z)}.
\]

It is worth mentioning that we use the ordinary generating functions of Stirling numbers instead of their exponential generating functions, which have been extensively used in \cite{Christoph2011}. Let
\[
F_1(n,k,z)= \frac{(z)_{\underline{n}}}{z^{k+1}} \quad \mbox{and} \quad
F_2(n,k,z) = \frac{z^k}{z^{n+1}(1-z)(1-2z)\cdots(1-k z)}.
\]
We see that both $F_1(n,k,z)$ and $F_2(n,k,z)$ are hypergeometric terms of $n$ and $k$. Let $N$ and $K$ be the shift operators with respect to $n$ and $k$, respectively. Denote the ring of linear difference operators with rational coefficients by
\[
\mathbb{K}(n,k)\langle N,K\rangle = \left\{ \sum_{i=0}^I \sum_{j=0}^J r_{i,j}(n,k) N^i K^j \colon I,J \in \mathbb{N},\ r_{i,j}(n,k) \in \mathbb{K}(n,k) \right\}.
\]
We see also that
\[
\Ann S_1(n,k) \subset \Ann F_1(n,k,z) \quad \mbox{and} \quad  \Ann S_2(n,k) \subset \Ann F_2(n,k,z),
\]
where
\[
\Ann f(n,k) = \{L \in \mathbb{K}(n,k)\langle N,K\rangle \colon L f(n,k) = 0 \}.
\]
Given a function $F(\n{n},k)$, we denote by $\tilde{F}(\n{n},k,z)$ the function obtained from $F(\n{n},k)$ by replacing $S_1(n,k)$ and $S_2(n,k)$ with $F_1(n,k,z)$ and $F_2(n,k,z)$, respectively. Suppose that there exist $Q\in \mathbb{K}(\textbf{n},k)\langle\textbf{N},K\rangle$ and $L\in \mathbb{K}(\textbf{n})\langle\textbf{N}\rangle$ such that
\[
L -(K-1)Q \in \Ann F(\textbf{n},k).
\]
Then we also have
\[
L -(K-1)Q \in \Ann \tilde{F}(\textbf{n},k,z),
\]
which leads to an equation of the form \eqref{tF-rec}.
The extended Zeilberger's algorithm will succeed in finding
such $L$ and $Q$. This fact indicates that the residue method always works as long as the existence of
such $L$ and $Q$ is guaranteed.

With the residue method, we can prove all identities on Stirling numbers appeared in \cite{Kauers2007}. Moreover, we can deal with sums involving products of Stirling numbers, typically are identities (6.24), (6.25), (6.28) and (6.29)
in \cite{GKP1994}. Here we only give two examples.

\begin{exam}
We have \cite[identity (6.24)]{GKP1994}
\begin{equation}\label{GKP6.24}
\sum_k S_1(k,m)S_2(n+1,k+1)=\binom{n}{m}.
\end{equation}
\end{exam}
\pf
Denote the left hand side by $L(n,m)$. We have
\[
L(n,m)=\res_x \res_y \sum_k \frac{x^{k+1}}{x^{n+2}\prod_{i=1}^{k+1}(1-ix)}\frac{(y)_{\underline{k}}}{y^{m+1}}.
\]
For the inner summand $F(n,m,k)$, Gosper's algorithm gives
\[
F(n,m,k)=G(n,m,k+1)-G(n,m,k),
\]
where
\[
G(n,m,k)=\frac{x^{k}}{x^{n+1}\prod_{i=1}^{k}(1-ix)}\frac{(y)_{\underline{k}}}{y^{m+1}}\frac{1}{1-x(1+y)}.
\]
Since the denominator contains $1-x(1+y)$ as a factor, we are unable to deduce a closed form of $\res_x \res_y G$. However, summing over $k$ from $0$ to $n$, we get
\[
L(n,m)= \res_x \res_y \frac{1}{(1-x(1+y))} \left(-\frac{1}{\prod_{i=1}^{n+1}(1-ix)}\frac{(y)_{\underline{n+1}}}{y^m}+\frac{1}{x^{n+1}y^{m+1}}\right).
\]
Notice that
\begin{align*}
\res_y \res_x \frac{1}{(1-x(1+y))}\frac{1}{\prod_{i=1}^{n+1}(1-ix)}\frac{(y)_{\underline{n+1}}}{y^m}=\res_y 0=0,
\end{align*}
and
\[
\res_y \res_x \frac{1}{(1-x(1+y))}\frac{1}{x^{n+1}y^{m+1}}=\res_y \frac{(1+y)^n}{y^{m+1}}=\binom{n}{m}.
\]
This completes the proof. \qed

\begin{exam}We have \cite[identity (6.28)]{GKP1994}
\begin{equation}\label{GKP6.28}
\sum_k \binom{n}{k}S_2(k,l)S_2(n-k,m)=\binom{l+m}{l}S_2(n,l+m).
\end{equation}
\end{exam}
\pf Denote the left hand side by $L(n,m,l)$. We have
\[
L(n,m,l)=\res_x \res_y \sum_k \binom{n}{k}\frac{x^l y^m}{x^{k+1}y^{n-k+1}\prod_{i=1}^{l}(1-ix)\prod_{j=1}^{m}(1-jy)}.
\]
For the inner summand $F(n,m,l,k)$, the extended Zeilberger's algorithm gives
\[
-F(n,m+1,l)-F(n,m,l+1)-(m+2+l)F(n,m+1,l+1)+F(n+1,m+1,l+1)=\Delta_k G(n,m,l,k),
\]
where
\[
G(n,m,l,k)=\binom{n}{k}\frac{-k x^l y^m}{(n+1-k)x^{k+1}y^{n-k+1}\prod_{i=1}^{l+1}(1-ix)\prod_{j=1}^{m+1}(1-jy)}.
\]
Summing over $k$ and applying the operators $\res_x$ and $\res_y$, we get a recurrence relation
\[
-L(n,m+1,l)-L(n,m,l+1)-(m+2+l)L(n,m+1,l+1)+L(n+1,m+1,l+1)=0.
\]
It is easy to check that the right hand side of \eqref{GKP6.28} satisfies the same recurrence relation.
Finally, the identities holds by checking the initial values
\[
L(0,m,l)=\delta_{0,l}\delta_{0,m}, \quad L(n,0,l)=S_2(n,l), \quad L(n,m,0)=S_2(n,m).
\]

\section{New identities}
In this section, we use two examples to illustrate how to discover new identities by the residue method.
In the first example, we generate new identities by introducing a new parameter in the original summand. While in the second example, we use Zeilberger's algorithm to construct new identities, as done by Chen and Sun \cite{ChenSun2009}.

We first consider the identity
\begin{equation}\label{AMM2011}
\sum_{k=0}^{n}(-1)^k \binom{2n}{n+k}S_1(n+k,k)=(2n-1)!!,
\end{equation}
which was proposed by Kauers and Sheng-Lang Ko as the American Mathematical Monthly Problem 11545. It was proved by F\"urst  \cite{Christoph2011} by the residue
representation
\[
S_1(n,k)=\frac{n!}{k!} \res_z \frac{\ln^k(1+z)}{z^{n+1}}.
\]
In fact, this identity can be generalized as follows.

\begin{theo}
Let $n$ and $m$ be nonnegative integers. Then we have
\begin{equation}\label{generalization-AMM2011}
\sum_{k=-m}^{n}(-1)^k \binom{2n}{n+k}S_1(n+k,k+m)=
\begin{cases}
(2n-1)!!, & \mbox{if $m=0$}, \\
0, & \mbox{if $m \ge 1$}.
\end{cases}
\end{equation}
\end{theo}
\pf Denote the left hand side of \eqref{generalization-AMM2011} by $L(n,m)$.
By the residue representation, we have
\[
L(n,m)=\res_z \sum_{k=-m}^{n} (-1)^k \binom{2n}{n+k} \frac{(z)_{\underline{n+k}}}{z^{m+k+1}}.
\]
The extended Zeilberger's algorithm gives the recurrence relation
\begin{multline}\label{AMMP-rec}
2(n+1)(2n+3)L(n,m)-\frac{(2n+3)(4n+3)}{2n+1}L(n+1,m)+L(n+2,m)\\
  +2(n+1)(2n+3)L(n+1,m+1)=0.
\end{multline}

We now prove that $L(n,n-r+1)=0$ for $n \ge r$ by induction on the non-negative integer $r$.
Since $S_1(n+k,n+1+k)=0$ for any integer $k$, we have $L(n,n+1)=0$, i.e., the assertion holds for $r=0$.
For $r=1$, we have
\[
L(n,n) = \sum_{k=-n}^n (-1)^k {2n \choose n+k} = 0, \quad n \ge 1.
\]
Now suppose that the assertion holds for $1 \le r \le r_0$ where $r_0 \ge 1$. The recurrence relation \eqref{AMMP-rec} implies that
\begin{multline*}
2(n-1)(2n-1)L(n-2,n-r_0)-\frac{(2n-1)(4n-5)}{2n-3}L(n-1,n-r_0)+L(n,n-r_0)\\
  +2(n-1)(2n-1)L(n-1,n-r_0+1) = 0.
\end{multline*}
By induction, we have
\[
L(n-2,n-r_0) = L(n-1,n-r_0) = L(n-1, n-r_0+1) = 0.
\]
Therefore, $L(n,n-r_0)=0$, which completes the induction. Notice that the assertion is equivalent to the statement $L(n,m)=0$ for any non-negative integers $n$ and $m \ge 1$.

For $m=0$, the recurrence relation \eqref{AMMP-rec} becomes
\[
2(n+1)(2n+3)L(n,0)-\frac{(2n+3)(4n+3)}{2n+1}L(n+1,0)+L(n+2,0)=0.
\]
It is easy to check that $(2n-1)!!$ satisfies this recurrence relation and coincides with the initial values $L(0,0)=L(1,0)=1$.
\qed

R. Sitgreaves \cite{Sitgreaves1970} found the
following identity (see also \cite{Egorychev1984}).
\begin{equation}\label{Sitgreaves}
\sum_{k=0}^{n}\binom{n+m}{k}(-1)^k S_2(n+m-k,n-k)=0, \quad m\geq 0, \quad n\geq m+1.
\end{equation}
From this result, we can establish the following theorem.
\begin{theo}
For nonnegative integers $n\geq m\geq 0$, we have
\begin{equation}\label{generalization-Sitgreaves}
\sum_{k=0}^{n}\binom{n+m+1}{k}(-1)^k S_2(n+m-k,n-k)=(-1)^{n+m}m!.
\end{equation}
\end{theo}
\pf Denoting the left hand side of \eqref{Sitgreaves} by $L(n,m)$, we have
\[
L(n,m)=\res_z \sum_{k=0}^{n} \binom{n+m}{k}(-1)^k \frac{z^{n-k}}{z^{n+m-k+1}\prod_{i=1}^{n-k}(1-iz)}.
\]
For the inner summand $F(n,m,k)$, the original Zeilberger's algorithm gives
\begin{multline*}
(n+m+1)F(n,m,k)-z(m+1)F(n,m+1,k)-zF(n,m+2,k)\\
=\Delta_k \frac{(n+m+1)kF(n,m,k)}{(n+m+1-k)(n+m+2-k)z}.
\end{multline*}
Summing over $k$ and applying the residue operator, we obtain
\begin{multline*}
(n+m+1)L(n,m)-(m+1)\sum_{k=0}^{n}\binom{n+m+1}{k}(-1)^k S_2(n+m-k,n-k)\\
-\sum_{k=0}^{n}\binom{n+m+2}{k}(-1)^k S_2(n+m+1-k,n-k)=0.
\end{multline*}
Denote the left hand side of \eqref{generalization-Sitgreaves} by $S(n,m)$.
Substituting $L(n,m)=0$ in the above identity, we deduce that
\[
(m+1)S(n,m)+S(n,m+1)=0.
\]
Thus we have
\[
S(n,m+1)=(-1)^{m+1}(m+1)!S(n,0).
\]
Note that
\[
S(n,0)=\sum_{k=0}^{n}\binom{n+1}{k}(-1)^k S_2(n-k,n-k)=\sum_{k=0}^{n}\binom{n+1}{k}(-1)^k =(-1)^n.
\]
we finally derive that
\[
S(n,m+1)=(-1)^{n+m+1}(m+1)!,
\]
as desired.\qed

\section{More combinatorial sequences}

It is readily seen that our approach is also appliable to many other combinatorial sequences as long as the
corresponding generating function is hypergeometric. More generally, the residue operator can be replaced by
any linear operator $L$. For example, a classical treatment for identities
involving harmonic numbers $H_n=\sum_{k=1}^n 1/k$ (see \cite{Paule-Schneider2003}) is to use the fact
\[
H_n=\delta \binom{n+x}{x},
\]
where $\delta f(x)=\frac{df(x)}{dx}|_{x=0}$.

Here we list several sequences which could be treated by this method.

\noindent \textbf{$q$-Stirling numbers}

A kind of $q$-analogue of Stirling numbers is given by \cite{Cigler1979,Johnson1994}
\begin{align*}
&S^{q}_1(n,k)=S^{q}_1(n-1,k-1)-[n-1]S^{q}_1(n-1,k),  &S^{q}_1(0,k)=\delta_{0,k},\\
&S^{q}_2(n,k)=S^{q}_2(n-1,k-1)+[k]S^{q}_2(n-1,k),  &S^{q}_2(0,k)=\delta_{0,k},
\end{align*}
where
\[
[n]=\frac{1-q^n}{1-q} = 1+q+\cdots+q^{n-1}.
\]
Their generating functions are
\[
\sum_{k=0}^{n}S^{q}_1(n,k) z^k =\prod_{k=0}^{n-1}(z-[k]), \quad
\sum_{n\geq k}S^{q}_2(n,k) z^n =\frac{z^k}{\prod\limits_{i=1}^{k}(1-[i]z)}.
\]
Thus we have
\[
S^{q}_1(n,k)=\res_z \frac{\prod\limits_{i=0}^{n-1}(z-[i])}{z^{k+1}}, \quad S^{q}_2(n,k)=\res_z \frac{z^{k-n-1}}{\prod\limits_{i=1}^{k}(1-[i]z)}.
\]
Note that we also have
\[
\Ann S^{q}_1(n,k) \subset \Ann \frac{\prod\limits_{i=0}^{n-1}(z-[i])}{z^{k+1}} \quad
\mbox{and} \quad  \Ann S^{q}_2(n,k) \subset \Ann \frac{z^{k-n-1}}{\prod\limits_{i=1}^{k}(1-[i]z)}.
\]

Using these representations and the $q$-analogue of the extended Zeilberger's algorithm, we can derive
recurrence relations for sums involving $q$-Stirling numbers.
For instance, let us consider the sum (see \cite{KauersSchneider})
\[
L(n,m)=\sum_{k=m}^{n}(-1)^{n-k}{k \brack m}_q S^{q}_1(n,k)q^{-k},
\]
where
\[
{k \brack m}_q = \frac{[k][k-1] \cdots [k-m+1]}{[m][m-1] \cdots [1]}
\]
is the $q$-binomial coefficients.
Our approach gives the recurrence relation
\[
L(n,m)+\frac{q(q^{m+1}-q^m+q^n-1)}{q-1}L(n,m+1)-q L(n+1,m+1)=0.
\]
Similarly, for the sum
\[
L(n,m)=\sum_{k=0}^{n}(-1)^{n-k}{n \brack k}_q S_{2}^{q}(k,m) q^{-k},
\]
we have
\begin{multline*}
(1-q^{n+1})L(n,m)+L(n+1,m)+\frac{(1-q^{m+1})(1-q^{n+1})}{1-q}L(n,m+1)\\
+\frac{(q^{m+1}-q^2+q-1)}{q-1}L(n+1,m+1)-qL(n+2,m+1)=0.
\end{multline*}

\noindent \textbf{Exponential functions}

Noting that
\[
k^n=\res_x \frac{1}{(1-kx)x^{n+1}},
\]
we can use the residue method to deal with sums involving $k^n$.
For example, consider the sum
\[
L(n,m)=\sum_k {m\choose k} k^n (-1)^{m-k}.
\]
Applying the extended Zeilberger's algorithm to the summand
\[
F(n,m,k)={m\choose k} \frac{(-1)^{m-k}}{(1-kx)x^{n+1}},
\]
we find that
\[
-(m+1)L(n,m)-(m+1)L(n,m+1)+L(n+1,m+1)=0.
\]
Since $m!S_2(n,m)$ satisfies the same recurrence relation and has the same
initial values, we finally derive that (see \cite[identity (6.19)]{GKP1994})
\[
\sum_k {m\choose k} k^n (-1)^{m-k} = m!S_2(n,m).
\]

\noindent \textbf{Bernoulli polynomials}

Identities involving Bernoulli and Euler numbers have been verified in \cite{ChenSun2009}. Here we only point out that we may also use the extended Zeilberger's
algorithm to derive differential equations satisfied by the sum. We take the Bernoulli polynomial
$B_n(x)$ as an example. Recall the generating function
\[
\sum_{n=0}^{\infty}B_n(x)\frac{z^k}{k!}=\frac{ze^{xz}}{e^{z}-1}.
\]
We have
\[
L(n,x,y)=\sum_{k}\binom{n}{k}y^{n-k}B_k(x) = \res_z \frac{1}{e^{z}-1}\sum_{k}\binom{n}{k}y^{n-k}\frac{e^{xz}k!}{z^k}.
\]
The extended Zeilberger's algorithm generates
\[
\frac{\partial}{\partial x} F(n+1,k,x,y)-(n+1)F(n,k,x,y)= \Delta_k \left( -\frac{(n+1)(n+1-k-ty)}{t(n+1-k)} F(n,k,x,y) \right).
\]
We thus have
\[
\frac{\partial}{\partial x} L(n+1,x,y) = (n+1) L(n,x,y).
\]
This relation together with the fact $L(n,0,y)=B_n(y)$ indicates that
\[
\sum_{k}\binom{n}{k}y^{n-k}B_k(x)=B_n(x+y).
\]

\section{Applicability of the residue method}

We have shown in Section 3 that for sums involving Stirling numbers, the residue method succeeds if the creative telescoping algorithm works whereas the converse is uncertain. In this section, we consider sums of the form
\[
\sum_k F(\n{n},k) a_k,
\]
where $F(\n{n},k)$ is a hypergeometric term and the generating function of $a_k$ is independent of $k$. By the residue method, we aim to find a finite set $S$ and $(k,z)$-free polynomial coefficients $\{ p_\bfal(\n{n}) \}_{\bfal \in S}$ such that
\[
\sum_{\bfal \in S} p_\bfal(\n{n}) \frac{F(\n{n}+\bfal, k)}{z^{k+1}} = \Delta_k G(\n{n},k,z).
\]
We will show that in most cases, the above equation holds only for $G(\n{n},k,z)=0$. In this case, we have
\[
\sum_{\bfal \in S} p_\bfal(\n{n}) F(\n{n}+\bfal, k) = 0,
\]
which is exactly the equation appears in Sister-Celine's method.

We first give a lemma on the $C$-finiteness of hypergeometric terms.
\begin{lem}\label{cfinite}
Let $f(k)$ be a hypergeometric term and
\begin{equation} \label{f-GP}
\frac{f(k+1)}{f(k)}=u\frac{A(k)}{B(k)}\frac{C(k+1)}{C(k)}
\end{equation}
be the GP-representation ((see \cite{PWZ1996} for the definition). If $f(k)$ is $C$-finite, then
\[
A(k) = B(k) = 1.
\]
\end{lem}
\pf Suppose that $f(k)$ is $C$-finite, this is, there exist constants $a_0,a_1,\ldots,a_d$, not all zeros, such that
\[
a_0 f(k)+a_1 f(k+1)+\cdots+a_d f(k+d)=0.
\]
Dividing $f(k)$ on both sides and substituting \eqref{f-GP}, we derive that
\begin{multline*}
a_0+a_1 u z\frac{A(k)}{B(k)}\frac{C(k+1)}{C(k)} + a_2 u^2 z\frac{A(k)A(k+1)}{B(k)B(k+1)}\frac{C(k+1)}{C(k)}+\cdots \\
+ a_d u^d \frac{\prod_{i=0}^{d-1}A(k+i)}{\prod_{i=0}^{d-1}B(k+i)}\frac{C(k+i)}{C(k)}=0.
\end{multline*}
Hence,
\[
\sum_{i=0}^{d}a_i u^i C(k+i)\prod_{j=0}^{i-1}A(k+j)\prod_{j=i}^{d-1}B(k+j)=0.
\]
Since $A(k)$ divides all the terms of the left hand side except the first one, it must also divides the first term. By the definition of GP-representation, $A(k)$ is co-prime to $C(k)$ and $B(k+j)$. We thus deduce that $A(k)=1$. With a similar discussion, we derive that $B(k) = 1$. \qed

Now we are ready to give the main theorem.
\begin{theo}\label{criterion}
Let $\n{n}=(n_1,\ldots,n_r)$ and $F(\n{n},k)$ be a hypergeometric term. Suppose that there are a finite set $S \subseteq \mathbb{N}^r$ and $(k,z)$-free polynomial coefficients $\{ p_\bfal(\n{n}) \}_{\bfal \in S}$ such that
\begin{equation}\label{ext_zeil_function}
\sum_{\bfal \in S} p_\bfal(\n{n}) \frac{F(\n{n}+\bfal, k)}{z^{k+1}} = \Delta_k G(\n{n},k,z).
\end{equation}
Let
\[
g(k) = \sum_{\bfal \in S} p_\bfal(\n{n}) F(\n{n}+\bfal, k),
\]
and
\[
\frac{g(k+1)}{g(k)} = u \frac{A(k)C(k+1)}{B(k)C(k)}
\]
be the GP-representation. Then we have $A(k)=B(k)=1$.
\end{theo}

\pf Since $F(\n{n},k)$ is hypergeometric, there exists a rational function $R(k,z)$ (since $\n{n}$ is irrelevant, we omit these variables) such that
\[
G(\n{n},k,z)= R(k,z) F(\n{n},k,z).
\]
Multiplying both sides of \eqref{ext_zeil_function} by
$z^{k+1}/F(\n{n},k)$, we see that
\[
h(k) = r(k) \frac{R(k+1,z)}{z}-R(k,z)
\]
is independent of $z$, where
\[
r(k) = \frac{F(\n{n},k+1)}{F(\n{n},k)}.
\]
Suppose that $R(k,z)=P(k,z)/Q(k,z)$, where $P(k,z)$ and $Q(k,z)$
are relatively prime polynomials in $k,z$. Then
\begin{equation}\label{pq-equation}
r(k)P(k+1,z)Q(k,z)-zP(k,z)Q(k+1,z)=h(k)zQ(k+1,z)Q(k,z).
\end{equation}
Noting that $r(k)$ and $h(k)$ are independent of $z$, by comparing the degrees in $z$ of both sides, we obtain $\deg_z P(k,z)=\deg_z Q(k,z)$.

We first prove that $z\nmid Q(k,z)$. Suppose on the contrary that there is a positive integer $m$ such
that $z^m \mid Q(k,z)$ but $z^{m+1}\nmid Q(k,z)$. By \eqref{pq-equation}, we see that $z^{m+1}\mid P(k+1,z)Q(k,z)$. Therefore, $z \mid P(k+1,z)$ and hence, $z\mid P(k,z)$. But this contradicts the condition that $P(k,z)$ and $Q(k,z)$ are relatively prime.

Then we show that $Q(k,z)$ is independent of $k$.
For any irreducible factor $p(k,z)$ of $Q(k,z)$, we deduce from \eqref{pq-equation} that $p(k,z)\mid z P(k,z)Q(k+1,z)$. Since
$z\nmid Q(k,z)$ and $P(k,z),Q(k,z)$ are relatively prime, we have $p(k,z)\mid Q(k+1,z)$, which implies $p(k-1,z)\mid Q(k,z)$. By iterating the above discussion, we get $p(k-i,z)\mid Q(k,z)$ for any nonnegative integer $i$. Therefore $p(k,z)$ must be independent of $k$. Since $p(k,z)$ is an arbitrary factor of $Q(k,z)$, we obtain that $Q(k,z)$ is independent of $k$.

From \eqref{pq-equation}, we see that $z\mid P(k+1,z)$. So we assume that
\[
P(k,z)=z\sum_{i=0}^{d}p_i(k)z^i, \quad Q(k,z)=\sum_{i=0}^{d+1}q_i z^i,
\]
where all $q_i$ are independent of $k$. Substituting these expressions into \eqref{pq-equation} and comparing the coefficient of each power of $z$, we find that
\begin{align}
r(k)p_0(k+1)&=q_0 h(k),\label{t1}\\[5pt]
r(k)p_1(k+1)-p_0(k)&=q_1 h(k)\label{t2},\\[5pt]
 &\vdots\\[5pt]
-p_d(k)&=q_{d+1} h(k).\label{td1}
\end{align}
By \eqref{t1}, we have $p_0(k+1)=q_0 h(k)/r(k)$. Substituting it into \eqref{t2}, we get
\[
p_1(k+2) = q_0 \frac{h(k+1)}{r(k)r(k+1)} + q_1 \frac{h(k+1)}{r(k)}.
\]
Continuing this discussion, we finally derive that
\[
q_0 h(k)+q_1 r(k)h(k+1)+q_2 r(k)r(k+1)h(k+2)+\cdots+q_{d+1}r(k)r(k+1)\cdots r(k+d) h(k+d+1)=0.
\]
Thus the hypergeometric term
\[
v(k)=h(k)\prod_{i=0}^{k-1}r(i)
\]
is $C$-finite.
Clearly,
\[
\frac{v(k+1)}{v(k)}= \frac{h(k+1)}{h(k)}r(k) = \frac{g(k+1)}{g(k)}.
\]
By Lemma~\ref{cfinite}, we deduce that $A(k)=B(k)=1$. \qed

For example, we consider the sum
\[
L(n,m) = \sum_{k=0}^n {n \choose k} B_{m+k},
\]
where $B_n$ is the Bernoulli number defined by
\[
\sum_{n=0}^{\infty}B_n\frac{z^n}{n!}=\frac{z}{e^{z}-1}.
\]
Rewrite the sum as
\[
L(n,m) = \sum_{k=m}^{n+m} {n \choose k-m} B_k.
\]
Denote the inner summand by $F(n,m,k)$. Take
\[
S=\{(i,j) \colon 0 \le i,j \le 1\},
\]
and denote
\[
g(k) = \sum_{(i,j) \in S} p_{i,j}(n,m) F(n+i,m+j, k).
\]
We find that
\[
\frac{g(k+1)}{g(k)}
= - \frac{k-m-n-2}{k+1-m} \frac{P(k+1)}{P(k)},
\]
where $P(k)$ is a certain polynomial in $k$. By Theorem~\eqref{ext_zeil_function}, we must have $g(k)=0$. There is a non-trivial solution
\[
p_{0,0}=p_{0,1}=-p_{1,0},
\]
which implies that
\[
L(n,m)+L(n,m+1)-L(n+1,m)=0.
\]
This coincides with the recurrence relation given by Chen and Sun \cite{ChenSun2009}, wherein all identities
involving only one Bernoulli number are of this case.

\noindent\textbf{Acknowledgments.} We thank Koutschan for his helpable discussion and suggestion. This work was supported by the PCSIRT project of the Ministry of
Education and the National Science Foundation of China.


\end{document}